\documentclass[10pt,reqno]{amsart}

\usepackage{amsaddr}
\usepackage{graphicx}
\usepackage{amssymb}
\usepackage{amsfonts}
\usepackage[]{amsmath}
\usepackage[]{epsfig}
\usepackage[]{float}
\usepackage{setspace}
\usepackage{tikz}
\usepackage{pgfplots}

\newtheorem{remark}{Remark}[section]

\numberwithin{equation}{section}
\newcommand{\R}{\mathbb R}

\newcommand{\intr}{\int_0^{+\infty}}

\usepackage[english, activeacute]{babel}
\usepackage{amsmath,amsthm,amsxtra}
\usepackage{epsfig}
\usepackage{amssymb}
\usepackage{latexsym}
\usepackage{amsfonts}
\usepackage{hyperref}
\usepackage{pdftricks}
\newtheorem{thm}{Theorem}[section]

\theoremstyle{remark}

\numberwithin{equation}{section}

\newcommand{\Z}{\mathbb{Z}}

\def\bm{\left( \begin{array}{cc}}
\def\endm{\end{array}\right)}

\newcommand{\be}{\begin{equation}}
\newcommand{\ee}{\end{equation}}
\newcommand{\ba}{\left(\begin{array}{c}}
\newcommand{\ea}{\end{array}\right)}
\newcommand{\bea}{\begin{eqnarray}}
\newcommand{\eea}{\end{eqnarray}}
\newcommand{\bee}{\begin{eqnarray*}}
\newcommand{\eee}{\end{eqnarray*}}
\newcommand{\ben}{\begin{enumerate}}
\newcommand{\een}{\end{enumerate}}

\begin{document}
\title[ Propagation of Regularity  of
Solutions to KDV on the Half-line]{On the Propagation of Regularity of
	Solutions to the KdV Equation on the positive Half-line}

\author{M\'arcio CAVALCANTE}
\address{\emph{Instituto de Matem\'atica, Universidade Federal de Alagoas,\\ Macei\'o-Brazil}}
\email{marcio.melo@im.ufal.br}
\author{Ailton C. Nascimento}
\address{\emph{Universidade Federal do Piauí,\\ Teresina-Brazil}}
\email{ailton.nascimento@ufpi.edu.br}


\begin{abstract}
We study special regularity properties of solutions to the initial-boundary value problem associated with the Korteweg-de Vries equations posed on the positive half-line. In particular, for initial data $u_0 \in H^{\frac{3}{4}^{+}}(\mathbb{R}^+)$  and boundary data $f\in H^{\frac32^+}(\R^+)$, where the restriction of $u_0$ to some subset of $(b,\infty)$ has an extra regularity for any $b>0$, we prove that the regularity of solutions $u$ moves with infinite speed to its left as time evolves until a certain time $T^*$. The existence of a stopping time  $T^{*}$ appears because of the effect of the boundary function $f$.  Also, as a consequence of our proof, we prove a gain in the regularity of the trace derivatives of the solutions for the Korteweg-de Vries on the half-line.\color{black}
\end{abstract}
\maketitle

\section{Introduction}

In this paper, we  study of the initial-boundary-value problem (IBVP) for the Korteweg-de Vries (KdV) equation posed on the half-line, namely
\begin{equation}\label{IBVP}
	\begin{cases}
		\partial_tu+\partial_x(\partial_x^2 u+u^2)=0  ,& (x,t)\in \mathbb R^+\times(0,T),\\
		u(x,0)=u_0(x),                                   & x\in \mathbb R^+,\\
		u(0,t)=f(t),                                     & t\in(0,T).
	\end{cases}
\end{equation}

\medskip

The IBVP \eqref{IBVP} is usually considered in the following setting: for $s\in\R$,
\begin{equation}\label{setting_A}
(u_0,f)\in H^{s}(\R^+)\times H^{\frac{s+1}{3}}(\R^+).
\end{equation}
This assumption is in some sense sharp because of the following localized smoothing effect for the linear evolution \cite{KPV}
\begin{align*}
&\|\psi(t) e^{-t\partial_x^3}\phi(x)\|_{C\big(\mathbb{R}_x;\; H^{(s+1)/3}(\mathbb{R}_t)\big)}\lesssim \|\phi\|_{H^s(\mathbb{R})},
\end{align*}
where $\psi(t)$ is a smooth cutoff function and $e^{-t\partial_x^3}$ denotes the linear homogeneous Airy group on $\mathbb{R}$. 

As pointed out by Zabuski \cite{Zabusky}, IBVP of the form \eqref{IBVP} may serve
as models for waves generated by a wave maker in a channel, or for waves approaching shallow water from deep water.  Here a mathematical issues connected to \eqref{IBVP} will be addressed; the propagation of regularity of solutions to the IBVP \eqref{IBVP} in the Sobolev spaces $H^s(\R^+).$ 

We begin with a review of existing theory which provides a setting in which to state precisely our results and put them into present day context.
The IBVP \eqref{IBVP} has been extensively studied in past few years, following the works by Ton \cite{Ton}, Bona and Winther \cite{BW,BW2}, Faminskii \cite{Fa,Fa2}, Bona, Sun and Zhang \cite{BSZ1,BSZ2}, Colliander and Kenig \cite{CK}, Holmer \cite{Holmer}, and Fokas \cite{Fokas1} (see also  \cite{Fokas2}). In particular, global well-posedness (GWP) does hold for data in $H^s(\R^+)$ for $s\geq 0$ and natural boundary conditions \eqref{setting_A}.

However, the behaviour of solutions to IBVP for the KdV equation has been considerably less studied than the corresponding initial value problem on complete real line. We now briefly comment on the results concerning orbital and asymptotic stability of solitons on the half-line. In \cite{CM} it was proved that the above-mentioned half-line solitons localized far from origin are orbitally stable in the energy space with respect to the flow of the IBVP \eqref{IBVP} for homogeneous boundary conditions. More recently, in the case of boundary conditions, in \cite{CM2} was obtained the asymptotic stability in the energy space, and were provided decay properties for all remaining regions, except the small soliton region. {For the KdV equation with a dominant surface tension, explicit solutions were obtained in \cite{Lennels} for appropriate initial and boundary conditions. While for the modified KdV equation on the half-line, a result of orbital stability for the Breathers was proved in \cite{Alejo}}.

{Finally, we point out that some recent results about various aspects associated with the behavior of solutions of the IBVPs for other nonlinear dispersive equations on the half-lines have been obtained in recent years. For the Schrödinger equation with nonlinear boundary conditions, Ackleh and Deng \cite{Ackleh} obtained blow-up solutions on the energy spaces and negative initial energy in some regimes depending on the structure of boundary conditions. In addition, blow-up for the 1D nonlinear Schrödinger equation with point nonlinearity was studied by Liu and Holmer \cite{Liu}, while Adami, Fukuizumi and Holmer \cite{Adami} considered the scattering problem of such a formulation. Furthermore, Kalantarov and Ozsari \cite{Kalantarov} proved that the solutions of the nonlinear Schrödinger equation  on the half-line, in some regimes, blow up in finite time by assuming the classical weighted condition on initial data. Also, Hayashi, Ogawa and Sato \cite{Sato} have considered the blow-up solutions in finite time for the nonlinear Schrödinger equation on the half-line with a nonlinear Neumann boundary condition, without the weight condition on the initial data. Finally, we point out that Chatziafratis, Ozawa and Tian \cite{Tian1} discovered a novel long-range instability phenomenon, which is a previously-unknown type for the inhomogeneous linear Schrödinger equation on the vacuum spacetime quarter-plane by developing the linear Fokas' unified transform method (see also the related recent works \cite{Tian3} and \cite{Tian2}). }

\color{black}

\subsection{Motivation and formulation of the problem}


Remember that for
the initial-value problem (IVP) posed in all $\R$  for the KdV-equation given by
\begin{equation}\label{pure}
\left\{\begin{array}{l}
\partial_{t} u+\partial_{x}^{3} u+u^{k} \partial_{x} u=0, \quad x, t \in \mathbb{R}, \quad k \in \mathbb{Z}^{+} \\
u(x, 0)=u_{0}(x)
\end{array}\right.
\end{equation}

The cases $k = 1$ and $k = 2$ in \eqref{pure} correspond to the KdV and modified KdV
(mKdV) equations, respectively.

This IVP has been extensively studied in the last years, we refer to the reader for instance
to the work of Killip and Vişan \cite{Killip}  and references therein. For a complete and detailed account, see also the monograph by Linares
and Ponce \cite{LP}.

Now, we comment on the results concerning the propagation of regularity for some nonlinear dispersive models. Initially, Isaza, Linares and Ponce \cite{Isaza} obtained for the first time the propagation of regularity for solution of the IVP \eqref{pure}. More precisely, they proved that the unidirectional dispersion of the k-generalized KdV equation \eqref{pure} produces the following propagation of regularity phenomena \cite{Isaza}: if for some $l\in\Z^{+}$ and $b\in \R$
\begin{equation}\label{Dat1}
	\|\,\partial_x^l u_0\|^2_{L^2((b,\infty))}<\infty,
\end{equation}
then for positive times the corresponding local solution $u=u(x,t)$ satisfies
\begin{equation}\label{Sol1}
	\|\,\partial_x^l u(\cdot,t)\|^2_{L^2((a,\infty))}<\infty \ \ \mbox{for every} \ a\in\R.
\end{equation}
This result tell us that the regularity $\eqref{Dat1}$ moves with infinite speed to its left as time evolves.
Interestingly, this enhanced regularity is known to propagate instantaneously throughout the solution of the flow, a phenomenon consistently observed in numerous nonlinear dispersive models. This phenomenon is now widely recognized as the {\it principle of propagation of regularity}. It was subsequently examined in the context of the Benjamin-Ono equation \cite{ILP} and the Kadomtsev-Petviashvili equation \cite{ILP1}. Similar results have been established for a broad class of equations; see, for example, \cite{AOR, LEV, LMP, LP, Arg, Arg1, MPS, Na1, Na, Na3}. The extension to the setting where the initial data possess fractional regularity was first explored by Kenig, Linares, Ponce, and Vega \cite{KLPV} in the framework of the KdV equation.

Moreover, recent contributions such as \cite{Arg2, Arg3, AO} have demonstrated analogous results for the $n$-dimensional Zakharov-Kuznetsov equation and its fractional counterparts, utilizing more advanced analytical techniques, including a novel class of pseudo-differential operators. Related findings were also obtained in \cite{Nascimento} for a family of nonlocal, nonhomogeneous, nonlinear dispersive models. For an up-to-date and thorough review of results concerning propagation of regularity, we refer the reader to Linares and Ponce \cite{LP2}, along with the references cited therein.

To the best of our knowledge, the present work constitutes the first result addressing the phenomenon of propagation of regularity for nonlinear dispersive equations posed on the half-line.

\color{black}
\subsection{Principal result} Before stating our results we need to define the class of solutions to the IBVP \eqref{IBVP} to which it applies. Thus, we shall rely on the following well-posedness result
which is a consequence of the arguments deduced in \cite{BSZ2} and \cite{Fa}.
\begin{thm}\label{katot}
For $s>\frac{3}{4}$, if $u_0 \in H^{s}\left(\R^{+}\right)$and $f \in H_{l o c}^{\frac{s+1}{3}}\left(\R^{+}\right)$ satisfy certain compatibility conditions at $(x,t)=(0,0)$, then the IBVP \eqref{IBVP} admits a unique solution
$$
u \in C\left(0, T ; H^{s}\left(\R^{+}\right)\right) \cap L^{2}\left(0, T ; H_{l o c}^{s+1}\left(\R^{+}\right)\right),
$$
which satisfies the additional properties

\begin{equation}\label{reg-trace}
\partial_x^ju \in C(\R^+; H^{\frac{s-j+1}{3}}((0,T))),\ j\in \{0,1,2,3,...\};
\end{equation}
\begin{equation}\label{kato}\left(\sup _{0<x<+\infty} \int_{0}^{T}\left|\partial_{x}^{s+1} u(x, t)\right|^{2} d t\right)^{\frac{1}{2}} \leq C\left(\|u_0\|_{H^{s}\left(\R^{+}\right)}+\|f\|_{H^{\frac{s+1}{3}}(0, T)}\right);
\end{equation}
\begin{equation}\label{stri}
\left(\int_{0}^{T} \sup _{0<x<+\infty}\left|\partial_{x} u(\cdot, t)\right|^{4} d t\right)^{\frac{1}{4}} \leq C\left(\|u_0\|_{H^{s}(\R^+)}+\|f\|_{H^{\frac{s+1}{3}}(0, T)}\right)
\end{equation}
and
\begin{equation}
\left(\int_{0}^{+\infty} \sup _{0 \leq t \leq T}|u(x, t)|^{2} d x\right)^{\frac{1}{2}} \leq C\left(\|u_0\|_{H^{s}(\R^+)}+\|f\|_{H^{\frac{s+1}{3}}(0, T)}\right),
\end{equation}
where the constants depend only on $s$ and $T$.
\end{thm}

\begin{remark}
	The property \eqref{kato}, now known as Kato's smoothing effect, estimulated an extensive investigation of various smoothing properties associated with solving the KdV-equation and other dispersive wave equations. In fact, \eqref{kato} is crucial to obtain propagation of the regularity result for nonlinear dispersive equations in general (see the discussion in \cite{LPSmith}).
	\end{remark}

The following interesting questions arise naturally in this situation:

\begin{itemize}
\item Is it possible to get results of propagation of the regularity for the solutions of the IBVP on the level obtained in the context of IVP by Isaza, Linares and Ponce \cite{Isaza} in the context of half-lines?

\item What is the influence of the known boundary functions and the unknown traces $\partial_x^ku(0,t)$ to get this result?

\item Is it possible to get propagation for any time of the existence?

\item With the assumption of regularity on the portion of the half-line for the initial data, is there any gain of the regularity of the trace  functions if it compared with \eqref{reg-trace}?

\end{itemize}

In this work we investigate these questions. Our result is concerned with the propagation of regularity in the right hand side of the data for positive times. It affirms that this regularity moves with infinite
speed to its left as time evolves.

\begin{thm}\label{thm1}
 If $u_{0} \in H^{3 / 4^{+}}(\mathbb{R}^+)$ and $f\in H^{\frac{3}{2}^+}(\R^+)$ for some $l \in \mathbb{Z}^{+}, l \geq 1$ and $x_{0} \in \mathbb{R}^+$
	
	\begin{equation}\label{hipl}
	\left\|\partial_{x}^{l} u_{0}\right\|_{L^{2}\left(\left(x_{0}, \infty\right)\right)}^{2}=\int_{x_{0}}^{\infty}\left|\partial_{x}^{l} u_{0}(x)\right|^{2} d x<\infty,
	\end{equation}
	then  for any $v \geq 0, \epsilon>0$ and $R>0$ hold the following statements for the solutions $u(t)$ on $[0,T]$ of the IBVP \eqref{IBVP} given by Theorem \ref{katot}:
	\begin{itemize}
\item[(a)] (Propagation of regularity) There exists $C>0$ such that
	$$
	\sup _{0 \leq t \leq T^*} \int_{\max\{ x_{0}+\epsilon-v t,0\}}^{\infty}\left(\partial_{x}^{j} u\right)^{2}(x, t) d x<C,
	$$
	where $T^*=T$ if $j=1$ and $T^*=\min\left\{T,\frac{x_0+\epsilon}{v}\right\}$; $j=2, \ldots, l$.
	In particular, for all $t \in(0, T^*]$, the restriction of $u(\cdot, t)$ to any interval $\left(x_{0}, \infty\right)$ belongs to $H^{l}\left(\left(x_{0}, \infty\right)\right)$.
\item[(b)] (Kato smoothing effect)	For any $v\geq 0$, $\epsilon>0$, $R>0$ and $T^*=T$ if $j=1$ and $T^*=\min\left\{T,\frac{x_0+\epsilon}{v}\right\}$; $j=2, \ldots, l$.
	$$
	\int_{0}^{T^*} \int_{\max\{ x_{0}+\epsilon-v t,0\}}^{x_{0}+R-v t}\left(\partial_{x}^{l+1} u\right)^{2}(x, t) d x d t<C.
	$$
	
	\item[(c)] (Trace estimate) Additionally, for  $T^*=\min\left\{T,\frac{x_0+\epsilon}{v}\right\}$, we have the following estimate for the trace: 
	 \be
	 \int_{\frac{b+x_0}{v}}^{T^*}(\partial_{x}^{2}u(0,t))^2d t\leq C.
	 \ee
	
		\end{itemize}
		(The constant $C$ depends on the numbers  $\left\|u_{0}\right\|_{H^{{3/4}^{+}}(\R^+)}$ $\left\|\partial_{x}^{l} u_{0}\right\|_{L^{2}\left(\left(x_{0}, \infty\right)\right)}$, $\|f\|_{H^{\frac{3}{2}^+}(\R^+)},$ $l,$ $v$, $\epsilon$ and  $T$.)
\end{thm}

The argument of the proof initially follows the ideas of \cite{Isaza}, but we need to take care with the trace terms which in general are not controllable on the low regularity regime. The main ingredients are the use of some Kato smoothing effects combined with a interpolation inequality, which is essential to control some problematic traces terms.
 \begin{figure}[htp]\label{figure}
	\centering 
	\begin{tikzpicture}[scale=2.2]
		\draw[->] (0,0)--(3,0) node[below] {$\boldsymbol{x}$ };
		
		\draw[color=blue,very thick, <-] (2,-0.25)--(3,-0.25);
		
		\node at (2.5,-0.33) {direction of propagation };
		
		\draw[->] (0,0)--(0,0.5) node[right] {$\boldsymbol{u(x,t)}$ };
		\draw[->] (0,0)--(-1,-1) node[right] {$\boldsymbol{t}$ };
		\draw[ultra thick][-](-0.3,-0.5)--(3,-0.5);
		\node at (-0.22,-0.6){$x_0+\epsilon-vt$};
		\node at (1,-0.4){$u(t)\in H^{1}$};
		
		\draw[ultra thick](1.05,0)--(2.97,0);
		
		\draw[ultra thick][-](-0.9,-0.9)--(3,-0.9);

		\node at (1.05,-0.08){$x_0$};
		\node at (1.7,0.1){$u_0\in H^1$};
		
	\end{tikzpicture}
	\caption{{\small The filled lines denote the region where regularity  occurs. In the case $l=1$, depending of the time of the existence the propagation can be reached.}}
	\label{Figura-I}
\end{figure}
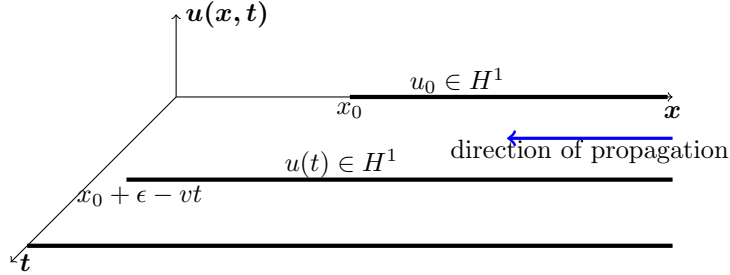

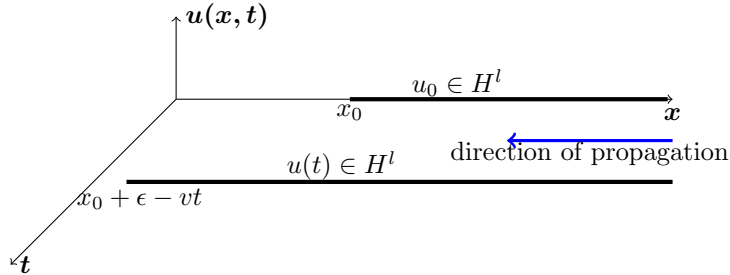
\begin{figure}[htp]\label{figure2}
	\centering 
	\begin{tikzpicture}[scale=2.2]
		\draw[->] (0,0)--(3,0) node[below] {$\boldsymbol{x}$ };
		
		\draw[color=blue,very thick, <-] (2,-0.25)--(3,-0.25);
		
		\node at (2.5,-0.33) {direction of propagation };
		
		\draw[->] (0,0)--(0,0.5) node[right] {$\boldsymbol{u(x,t)}$ };
		\draw[->] (0,0)--(-1,-1) node[right] {$\boldsymbol{t}$ };
		\draw[ultra thick][-](-0.3,-0.5)--(3,-0.5);
		\node at (-0.22,-0.6){$x_0+\epsilon-vt$};
		\node at (1,-0.4){$u(t)\in H^{l}$};
		
		\draw[ultra thick](1.05,0)--(2.97,0);

		\node at (1.05,-0.08){$x_0$};
		\node at (1.7,0.1){$u_0\in H^l$};
		
	\end{tikzpicture}
	\caption{{\small  In the case $l\geq 2$ the propagation cannot reach the boundary.}}
	\label{Figura-I}
\end{figure}

\begin{remark} (About the bound $\frac{x_0+\epsilon}{v}$) Unlike the KdV equation in the context of the entire line \cite{Isaza}, our propagation results have a new bound due to higher-order derivative trace functions. These functions cannot be estimated with the assumption of regularity beyond $\frac{3}{4}$. Therefore, this additional hypothesis regarding time prevents the presence of these unfavorable trace terms. 
	\end{remark}

\begin{remark}
	The propagation result on the case $l=1$ is more efficient since here we do not need to give a addidtional bound for the time. Also, in this situation the part (c) of Theorem \ref{thm1} gives a important gain of the regularity for the trace function $\partial_{x}^2u(0,t)$, wich give a gain of regularity for this function. More precisely the result on the literature gives $\partial_x^2u(0,t)$ on the space $H^{-\frac{1}{12}^+}(\R^+)$, while our results prove that   $\partial_{x}^2u(0,t)\in L^2(\frac{b+x_0}{v},\infty)$.
\end{remark}

\section{Preliminaries}

Let $\eta(\theta), \theta \in \mathbb{R}$, be a certain cut-off function, namely, $\eta \in C^{\infty}(\mathbb{R})$, $\eta \geq 0, \eta^{\prime} \geq 0, \eta(\theta)=0$ for $\theta \leq 0, \eta(\theta)=1$ for $\theta \geq 1, \eta(\theta)+\eta(1-\theta) \equiv 1$, $\eta^{\prime}(\theta)>0$ for $\theta \in(0,1)$.

Now, let $\rho(x) \equiv 1+\eta((x+2) / 3)$.

We will use the following auxiliary functions $\chi_{ \epsilon, b}(x)$ for $\epsilon>0$ and $b \geq 5 \epsilon$ such that
\be\label{chi}
\begin{gathered}
\chi_{ \epsilon, b} \in C^{\infty}(\mathbb{R}^+ ), \quad \chi_{ \epsilon, b}^{\prime} \geq 0, \\
\chi_{ \epsilon, b}(x)= \begin{cases}0, &0< x \leq \epsilon, \\
1, & x \geq b.\end{cases}
\end{gathered}
\ee

Note that the function $\chi_{\epsilon,b}'$ has compact support and it is bounded (see \cite{Isaza} for the properties of these weighted functions).

For $s\geq 0$ we say that $\phi \in H^s(\mathbb{R}^+)$ if there exists $\tilde{\phi}\in H^s(\mathbb{R})$ such that 
$\phi=\tilde{\phi}|_{\R+}$.  In this case we set $\|\phi\|_{H^s(\mathbb{R}^+)}:=\inf\limits_{\tilde{\phi}}\|\tilde{\phi}\|_{H^{s}(\mathbb{R})}$. 

\section{Start of the Proof: Case $l=1$ of the principal result}

Formally, take partial derivative with respect to $x$ of the equation in \eqref{IBVP} and multiply by $\partial_{x} u \chi_{\epsilon, b}(x+v t-x_0)$ to obtain after integration by parts the identity
\begin{equation}
\begin{aligned}
	&\frac{1}{2} \frac{d}{d t} \intr\left(\partial_{x} u\right)^{2}(x, t) \chi_{\epsilon,b}(x+v t-x_0) d x-\underbrace{v	 \intr\left(\partial_{x} u\right)^{2}(x, t) \chi_{\epsilon,b}^{\prime}(x+v t-x_0) d x}_{I_1} \\
	&\quad+\frac{3}{2} \intr\left(\partial_{x}^{2} u\right)^{2}(x, t) \chi_{\epsilon,b}^{\prime}(x+v t-x_0) d x-\underbrace{\frac{1}{2} \intr\left(\partial_{x} u\right)^{2}(x, t)\chi_{\epsilon,b}^{\prime \prime \prime}(x+v t-x_0) d x}_{I_{2}} \\
	&\quad+\underbrace{\intr \partial_{x}\left(u \partial_{x} u\right) \partial_{x} u(x, t)\chi_{\epsilon,b}(x+v t-x_0) d x}_{I_{3}}=\tau_1(t),
\end{aligned}
\end{equation}
where 
\begin{equation} 
\begin{split}
\tau_1(t)&= (\partial_xu(0,t))^2\chi_{\epsilon,b}''(vt-x_0) {-2\partial_{x}u(0,t)\chi_{\epsilon,b}'\partial_{x}^2u(0,t)}{-\frac12(\partial_{x}^2u(0,t))^2\chi_{\epsilon,b}(vt-x_0)}\\
&\quad\quad+\partial_{x}u(0,t)\chi_{\epsilon,b}(vt-x_0) \partial_x^3u(0,t) +\chi_{\epsilon,b}(vt-x_0) u(0,t) (\partial_xu)^2(0,t)\\&\quad\quad +\partial_x(u^2)(0,t)
\partial_xu(0,t) \chi_{\epsilon,b}(vt-x_0).\end{split}
\end{equation}

Now, by using the equation in \eqref{IBVP} we have that 
\begin{equation}\label{tec1}
\partial_x^3u(0,t)=-f'(t)-2f(t)\partial_x u(0,t),
\end{equation}
it is suffices to control the term that contains three derivatives. {In this step we need to assume a more regular boundary condition, more precisely we assume $f \in H^{{\frac32}^+}(\R^+)$, instead of the natural assumption $f\in H^{{\frac{7}{12}}^+}(\R)$.}

{The estimates of $I_1$ and $I_2$ follow from the arguments in  \cite{Isaza}.} In addition, for the $I_3$ term we have 
\begin{align*}
	I_3 &= \intr (\partial_xu)^3 \chi_{\epsilon,b}(x+v t-x_0)\,dx + \intr u \partial_x^2 u \partial_xu \chi_{\epsilon,b}(x+v t-x_0) \,dx \\
	&=  \frac12\intr (\partial_xu)^3 \chi_{\epsilon,b}\,dx-\frac12 \intr u(\partial_xu)^2\chi_{\epsilon,b}'dx-\frac12 u(0,t)(\partial_xu(0,t))^2\chi_{\epsilon,b}(vt-x_0)\\
	&:= I_{3_1} + I_{3_2} -\frac12 f(t)(\partial_xu(0,t))^2\chi_{\epsilon,b}(vt-x_0).
\end{align*}
These first two terms above can be controlled in the following way
\begin{equation}
	|I_{3_1}| \leq  \|\partial_xu \|_{L^{\infty}(\mathbb R^+)} \intr (\partial_x u)^2 \chi_{\epsilon,b}(x+vt-x_0) \,dx 
\end{equation}
and
\begin{equation}
	|I_{3_2}| \leq \|u(t)\|_{L^{\infty}(\mathbb R^+)} \intr (\partial_x u)^2 \chi_{\epsilon,b}'(x+vt-x_0) \,dx. \tag{3.6}
\end{equation}
Hence, by Sobolev embedding and (3.2) after integrating in the time interval $[0,T]$ one gets
\begin{equation}
	\int_0^T |I_{3_2}| \,dt \leq \sup_t \|u(t)\|_{H^{3/4 +}(\R^+)} \int_0^T \intr (\partial_x u)^2 \chi_{\epsilon,b}'(x+vt-x_0) \,dx dt \leq C_0. \tag{3.7}
\end{equation}
It follows that
\begin{equation}
\begin{split}
\sup _{0 \leq t \leq T*} &\intr\left(\partial_{x} u\right)^{2} \chi_{ \epsilon, b}(x+v t-x_0) d x+\int_{0}^{T^*} \intr\left(\partial_{x}^{2} u\right)^{2} \chi_{\epsilon, b}^{\prime}(x+v t-x_0) d x d t \\
&\leq C+\int_{0}^{T^*}\tau_1(t)dt,
\end{split}
\end{equation}
with $c=c(\epsilon ; b ; v;T^*)>0$ for any $\epsilon>0, b \geq 5 \epsilon, v>0$.

\medskip
Now, we treat each term of the trace $\tau_1(t)$. By using the regularity of the traces given in \eqref{reg-trace} we have that $u(0,t)\in H^{\frac{3}{2}^+}(0,T),\ \partial_x u(0,t)\in H^{\frac{1}{4}^+}(0,T)$ with controls the terms of $\tau_1$ which contain this terms.

Also, by using Young inequality we have for $\delta$ small enough
\begin{equation}
{-2\partial_{x}u(0,t)\chi_{\epsilon,b}'(vt-x_0)\partial_{x}^2u(0,t)\leq \delta (\chi_{\epsilon,b}'(vt-x_0)\partial_{x}^2u(0,t))^2+\frac{4}{\delta}(\partial_xu(0,t))^2}.
\end{equation}

Then,  we obtain the improved inequality
\begin{equation*}
	\begin{split}
\sup _{0 \leq t \leq T*} &\intr\left(\partial_{x} u\right)^{2} \chi_{ \epsilon, b}(x+v t-x_0) d x+\int_{0}^{T^*} \intr\left(\partial_{x}^{2} u\right)^{2} \chi_{\epsilon, b}^{\prime}(x+v t-x_0) d x d t \\
& {+\int_0^{T^*}\frac12(\partial_{x}^2u(0,t))^2(\chi_{\epsilon,b}(vt-x_0) -\delta(\chi_{\epsilon,b}')^2(vt-x_0))d t}\leq C+\int_{0}^{T^*}\widetilde \tau_1(t)dt,
\end{split}
\end{equation*}
where,
\begin{align}
	\widetilde \tau_1(t)&= (\partial_xu(0,t))^2\chi_{\epsilon,b}''(vt-x_0) +\partial_{x}u\chi_{\epsilon,b}(vt-x_0) \partial_x^3u(0,t) \nonumber\\
	&\quad \quad +\chi_{\epsilon,b}(vt-x_0) u(0,t) (\partial_xu)^2(0,t)+\partial_x(u^2)(0,t)
	\partial_xu (0,t)\chi_{\epsilon,b}(vt-x_0).
\end{align}

Now, by using  the assumption $u(0,t)=f(t)\in H^{\frac{7}{12}^+}(0,T),$ and $ \partial_x u(0,t)\in H^{\frac{1}{4}^+}(0,T)$ (from \eqref{reg-trace} ), we have that
\begin{equation}
\int_{0}^{T^*}\widetilde \tau_1(t)dt\leq C.
\end{equation}

Thus, we have proved statements (a) and (b) of Theorem \ref{thm1}, on the case $l=1$.
 
Finally, by using the properties of the functions $\chi_{\epsilon,b}$ we obtain the following localized smoothing effect
\begin{align}
	&\int_{\frac{b+x_0}{v}}^{T^*}\left (\frac12-\delta k\right)(\partial_{x}^2u(0,t))^2)d t \nonumber\\
	&\quad \quad \quad\leq \int_{\frac{b+x_0}{v}}^{T^*}\frac12(\partial_{x}^2u(0,t))^2(\chi_{\epsilon,b}(vt-x_0) -\delta(\chi'_{\epsilon,b})^2(vt-x_0))d t\leq C.
\end{align}

It follows that  
\begin{align}
	\int_{\frac{b+x_0}{v}}^{T^*}(\partial_{x}^2u(0,t))^2d t\leq C.
\end{align}

 which proves part (c) of Theorem \ref{thm1} on the case $l=1$.
\medskip

\textbf{Case 2: $l=2$}

By taking two derivatives in equation \eqref{IBVP} and multiplying by $\partial_x^2u \chi_{\epsilon,b}(x+vt-x_0)$ we obtain by making some integration by parts 

\begin{equation}\label{cl2}
\begin{aligned}
&\frac{1}{2} \frac{d}{d t} \intr\left(\partial_{x}^{2} u\right)^{2}(x, t) \chi_{\epsilon,b}(x+v t-x_0) d x-\underbrace{v \intr\left(\partial_{x}^{2} u\right)^{2}(x, t) \chi_{\epsilon,b}^{\prime}(x+v t-x_0) d x}_{I_{1}} \\
&\quad+\frac{3}{2} \intr\left(\partial_{x}^{3} u\right)^{2}(x, t) \chi_{\epsilon,b}^{\prime}(x+v t-x_0) d x-\underbrace{\frac{1}{2} \intr\left(\partial_{x}^{2} u\right)^{2}(x, t) \chi_{\epsilon,b}^{\prime \prime \prime}(x+v t-x_0) d x}_{I_{2}} \\
&\quad+\underbrace{\intr \partial_{x}^{2}\left(u \partial_{x} u\right) \partial_{x}^{2} u(x, t) \chi_{\epsilon,b}(x+v t-x_0) d x}_{I_{3}}+{(\partial_x^3u(0,t))^2\chi_{\epsilon,b}(v t-x_0)}=\tau_2(t),
\end{aligned}
\end{equation}
where
\begin{equation}
\begin{split}
\tau_2(\tau)&={-\partial_x^2u(0,t)\partial_x^4u(0,t) \chi_{\epsilon,b}(vt-x_0)+}\partial_x^2u(0,t)\partial_x^3u(0,t)\chi_{\epsilon,b}'(vt-x_0)\\
&\quad\quad-\frac{1}{2}(\partial_x^2u(0,t))^2\chi_{\epsilon,b}''(vt-x_0)+\widetilde \tau_2(t),
\end{split}
\end{equation}
where the function $\widetilde \tau_2(t)$ depends of  $\chi$  and of the traces functions $\partial_x^k u(0,t)$ for $k\in \{0,1,2,3\}$ and it was treated on the previous case.

Here we have a problematic term which contains the trace of fourth-order derivatives of $u$. Unlike the previous case, the equation does not help to control this problem, then we need to impose a bound for the time $t$. In this context, we will assume  $t< \frac{(\epsilon+x_0)}{v}$  which implies $vt-x_0<\epsilon$. With this assumption the problematic term ${-\partial_x^2u(0,t)\partial_x^4u(0,t) \chi_{\epsilon,b}'(vt-x_0)}$ disappears. This implies that in this case the regularity does not reach the boundary. The estimates of $I_1$ and $I_2$ follow from the arguments in  \cite{Isaza}.

Regarding the term $I_3$ on the LHS of \eqref{cl2}, we observe that
\begin{align}\nonumber
	I_3= & \intr \partial_x^{2}\left(u \partial_x u\right) \partial_x^{2} u \chi_{\epsilon,b}(x+v t-x_0) d x \\\label{l2at}
	= & \frac{5}{2}\intr \partial_x u\left(\partial_x^{2} u\right)^2 \chi_{\epsilon,b}(x+v t-x_0) d x- \frac{1}{2}\intr u \left(\partial_x^{2} u\right)^2  \chi_{\epsilon,b}^{\prime}(x+v t-x_0) d x \\\nonumber
	& -u(0,t)\left(\partial_x^{2} u(0,t)\right)^2  \chi_{\epsilon,b}(v t-x_0)\\\nonumber
:=&	I_{31}+I_{32}-\tau(t).
\end{align}

The $I_{31}$ term can be estimated as follows,

$$|I_{31}| \leq c  ||\partial_x u(t)||_{L^{\infty}(\mathbb R^+)} \intr (\partial_x^2u)^2 \chi_{\epsilon,b}(x+vt-x_0) dx,$$
where the last integral is the quantity to be estimated.

For the second term we have
$$\int_{0}^{T}|I_{32}(t)|dt \le \sup_{0 \le t \le T} ||u(t)||_{L^{\infty}(\mathbb R^+)} \int_{0}^{T} \intr (\partial_x^2 u)^2 \chi_{\epsilon,b}^{\prime}(x+vt-x_0)dxdt \le c.$$

Inserting the above information in \eqref{cl2} and using Gronwall's inequality 
\begin{equation*}
	\sup_{0 < t \le T} \intr (\partial_x u)^2\chi_{\epsilon,b} (x+vt-x_0) dx + \int_0^T \intr (\partial_x^2 u)^2 \chi_{\epsilon,b}'(x+vt-x_0) dx dt \le c_0,
\end{equation*}
where we have used that 
$$\int_0^T|\tau(t)|dt\leq c,$$  because of the previous step.

Finally, we have the following additional gain of regularity for the trace of the third derivative
\begin{equation}
	\int_{\frac{b+x_0}{v}}^{T^*}{(\partial_x^3u(0,t))^2dt}\leq C.
\end{equation}

Before to treat the general case, we illustrate the control of the nonlinear term for the case $l=3$. More precisely, we have the following estimate
\begin{align}\nonumber
	I_4:= & \int_0^{+\infty} \partial_x^3\left(u \partial_x u\right) \partial_x^3 u \chi_{\epsilon,b}(x+vt-x_0) dx \\\nonumber
	= & 4 \int_0^{+\infty} \partial_x u\left(\partial_x^3 u\right)^2 \chi_{\epsilon,b} dx+\int_0^{+\infty} u \partial_x^4 u \partial_x^3 u \chi_{\epsilon,b} dx  +3 \int_0^{+\infty} (\partial_x^2 u)^2 \partial_x^3 u \chi_{\epsilon,b} dx \\\nonumber
	=&\frac{7}{2} \int_0^{+\infty} \partial_x u\left(\partial_x^3 u\right)^2 \chi_{\epsilon,b} dx-\frac{1}{2} \int_0^{+\infty} u\left(\partial_x^3 u\right)^2 \chi'_{\epsilon,b}dx-\frac12\left.\left(\partial_x^3 u\right)^2(0,t) u \chi_{\epsilon,b}(vt-x_0)\right.\\ 
	&-\frac{3}{2}\intr\left(\partial_x^2 u\right)^3 \chi'_{\epsilon,b}dx-\left.\left(\partial_x^2 u\right)^3 (0,t)\chi_{\epsilon,b}(vt-x_0).\right.
\end{align}

In order to control the antepenultimate term, we invoke the properties of the weighted function to get
\begin{equation}\label{l32e}
\begin{split}
	&\left|\intr\left(\partial_x^2 u\right)^3 \chi'_{\epsilon,b}(x+vt-x_0)dx\right| \\
	&\quad\leq\left\|{\partial}_x^2 u \chi_{ \epsilon, b}^{\prime}(\cdot+v t-x_0)\right\|_{L^\infty(\R^+)} \intr\left({\partial}_x^2 u\right)^2 \chi_{ \epsilon / 5, \epsilon}(x+v t-x_0) d x,
\end{split}
\end{equation}
where the last term is bounded, after integration in time, by the former case $l=2$. Next, we take care of the first term on the RHS of \eqref{l32e} as follows

\begin{align}\nonumber
&	\left\|\left(\partial_x^2 u\right)^2 \chi_{ \epsilon, b}^{\prime}(\cdot+v t-x_0)\right\|_{\infty}  \lesssim \intr\left|\partial_x\left(\left(\partial_x^2 u\right)^2 \chi_{ \epsilon, b}^{\prime}(x+v t-x_0)\right)\right| d x \\\nonumber
	& \qquad\lesssim c \intr\left|\partial_x^2 u \partial_x^3 u \chi_{ \epsilon, b}^{\prime}(x+v t-x_0)\right| d x +\intr\left|\left(\partial_x^2 u\right)^2 \chi_{ \epsilon, b}^{\prime \prime}(x+v t-x_0)\right| d x\\\nonumber
	&\qquad\lesssim  c \intr\left(\partial_x^2 u\right)^2 \chi_{ \epsilon, b}^{\prime}(x+v t-x_0) d x+c \intr\left(\partial_x^3 u\right)^2 \chi_{ \epsilon, b}^{\prime}(x+v t-x_0) d x \\\label{l2arg}
	& \qquad\qquad+c \intr\left(\partial_x^2 u\right)^2 \chi_{ \epsilon / 3, b+\epsilon}^{\prime}(x+v t-x_0) d x.	
\end{align}

 \section{End of the proof: Induction Argument}

 We shall prove the case $l+1$ assuming the case $l \geq 2$. More precisely, we assume:
 If $u_0$ satisfies (1.13) then (1.14) holds, i.e.
 \begin{equation}\label{HIl1}
 	\sup _{0 \leq t \leq T} \intr\left(\partial_x^j u\right)^2 \chi_{\epsilon, b}(x+v t-x_0) d x+\int_0^T \intr\left(\partial_x^{j+1} u\right)^2 \chi_{\epsilon, b}^{\prime}(x+v t-x_0) d x d t \leq C
 \end{equation}
 for $j=1,2, \ldots, l, l \geq 2$, for any $\epsilon>0, b \geq 5 \epsilon, v>0$.
 
 Now we have that
 \begin{equation}\label{Hipl1}
 	\left. u_0\right|_{(x_0, \infty)} \in H^{l+1}([x_0, \infty)).	
 \end{equation}
 Thus from the previous step $\eqref{HIl1}$ holds. Also formally we have for $\varepsilon>0, b \geq 5 \epsilon$, the identity
 \begin{align}\nonumber
 	\frac{1}{2} & \frac{d}{d t} \intr\left(\partial_x^{l+1} u\right)^2 \chi_{\epsilon, b}(x+v t-x_0) d x-\underbrace{v \intr\left(\partial_x^{l+1} u\right)^2 \chi_{\epsilon, b}^{\prime}(x+v t-x_0) d x}_{I_1} \\\nonumber
 	+ & \frac{3}{2} \intr\left(\partial_x^{l+2} u\right)^2 \chi_{\epsilon, b}^{\prime}(x+v t-x_0) d x+\underbrace{\left(\partial_x^{l+2} u(0, t)\right)^2 \chi_{\epsilon, b}(v t-x_0)}_{I_2}\\\nonumber&\quad-\underbrace{\frac{1}{2} \intr\left(\partial_x^{l+1} u\right)^2 \chi_{\epsilon, b}^{\prime \prime \prime}(x+v t-x_0)  d x}_{I_3} \\\label{lcasel1}
 	&\quad +\underbrace{\intr \partial_x^{l+1}\left(u {\partial}_x u\right) \partial_x^{l+1} u(x, t) \chi_{\epsilon, b}(x+v t-x_0)  d x}_{I_4}+\tau_{l+1}(t)=0,
 \end{align}
 where the trace term $\tau_{l+1}(t)$ is given as
 \begin{equation}\label{tracel1}
 	\begin{split}
 	\tau_{l+1}(t)=&-\partial_x^{l+3} u(0, t) \partial_x^{l+1} u(0, t) \chi_{\epsilon, b}(v t-x_0) \\&\quad +\partial_x^{l+2} u(0, t) \partial_x^{l+1} u(0, t) \chi_{\epsilon,b}^{\prime}\left(v t-x_0\right)-\left(\partial_x^{l+1} u\right)^2 \chi_{\epsilon, b}^{\prime \prime \prime}(v t-x_0).	
 \end{split}
 \end{equation}
 We observe that all the term on the LHS of $\eqref{lcasel1}$ can be treated as we did previously in the former cases.  We also note that the trace term $\eqref{tracel1}$ is of order $l+1$ and can be estimated employing the argument used in the analysis of the cases $l=1$ and $l=2$, taking advantage of the properties of weighted function $\chi_{\epsilon,b}$. However, we illustrate the analysis of the nonlinear term $I_4$. In fact, we write
 \begin{align}\nonumber
 	I_4= & \intr \partial_x^{l+1}\left(u \partial_x u\right) \partial_x^{l+1} u \chi_{\epsilon,b}(x+v t-x_0) d x \\\nonumber
 	= & c_0 \intr u \left(\partial_x^{l+1} u\right)^2  \chi_{\epsilon,b}^{\prime}(x+v t-x_0) d x \\\nonumber
 	& +c_1 \intr \partial_x u\left(\partial_x^{l+1} u\right)^2 \chi_{\epsilon,b}(x+v t-x_0) d x \\\nonumber
 	& +c_2 \intr \partial_x^2 u \partial_x^l u \partial_x^{l+1} u \chi_{\epsilon,b}(x+v t-x_0) d x \\\label{A4l1}
 	& -\left(\partial_x^{l+1} u(0,t)\right)^2 u(0,t) \chi_{\epsilon, b}(v t-x_0)+\sum_{k=3}^{l-1} c_k\intr \partial_x^{l+2-k} u \partial_x^{k} u \partial_x^{l+1} u \chi_{\epsilon,b}(x+v t-x_0) d x \\\nonumber
 	:= & I_{4,0}+I_{4,1}+I_{4,2}+ \tilde\tau_{l+1}(t)+\sum_{k=3}^{l-1} I_{4, k}.
 \end{align}
 
 Thus, a familiar argument yields
 \begin{equation}\label{a40}
 	\left|I_{4,0}(t)\right| \leq\|u(t)\|_{L^\infty(\R^+)} \intr\left(\partial_x^{l+1} u\right)^2 \chi_{\epsilon,b}^{\prime}(x+v t-x_0) d x 	
 \end{equation}
 which after integration in time is bounded by Sobolev embedding and the hypothesis \eqref{HIl1} $j=l$ provides the boundeness.
 
 Also,
 \begin{equation}\label{a41}
 	\left|I_{4,1}(t)\right| \leq\left\|\partial_x u(t)\right\|_{L^\infty(\R^+)} \intr\left(\partial_x^{l+1} u\right)^2 \chi_{\epsilon,b}(x+v t-x_0) d x .
 \end{equation}
 We employ the estimate \eqref{stri} to bound the first term on the RHS of $\eqref{a41}$, after the integration in time. The remaining term is the very quantity to be estimated. The term $I_{4,2}$ is bounded following the same argument we employed before in the analysis \eqref{l32e}-\eqref{l2arg}.
 
 Finally, we see that
 \begin{align}\nonumber
 	\left|I_{4, k}(t)\right|= & \left|\intr \partial_x^{l+2-k} u \partial_x^{k} u \partial_x^{l+1} u \chi_{\epsilon,b}(x+v t-x_0) d x\right| \\\nonumber
 	\leq & \frac{1}{2} \intr\left(\partial_x^{l+2-k} u \partial_x^{k} u\right)^2 \chi_{\epsilon,b}(x+v t-x_0) d x \\\label{a4k}
 	& +\frac{1}{2} \intr\left(\partial_x^{l+1} u\right)^2 \chi_{\epsilon,b}(x+v t-x_0) d x 
 \end{align}
 where the last term is the quantity to be estimated. We notice that $k$, $l+2-k \leq l-1$ so we obtain
\begin{align}\label{fina4k}
	&\left|\intr\left(\partial_x^{l+2-k} u \partial_x^{k} u\right)^2 \chi_{\epsilon,b}(x+v t-x_0) d x\right| \leq \nonumber\\
	&\qquad\left\|\left(\partial_x^{k} u\right)^2 \chi_{ \epsilon / 5, \epsilon}(\cdot+v t-x_0)\right\|_{L^\infty(\mathbb R^+)} \intr\left(\partial_x^{l+2-k} u\right)^2 \chi_{\epsilon,b}(x+v t-x_0) d x.
\end{align}

 The first term on the RHS of \eqref{fina4k} can be treated similarly as we did in the analysis of the $I_{4,2}$ term in \eqref{A4l1}, taking into the account the fact that $k\leq l-1$. The last term is bounded by the hypothesis \eqref{HIl1} since $l+2-k \leq l-1$. This concludes the induction argument.
 
 \section*{Acknowledgments}
 M. Cavalcante wishes to thank the support of CNPq, 
 grant \# 310271/
 2021-5, and the support of the Funda\c c\~ao de Amparo \`a Pesquisa do Estado de Alagoas
 - FAPEAL, Brazil, grant \# E:60030.0000000161/2022. The authors thank the CAPES/Cofecub, grant \# 88887.879175/2023-00.

COI: All authors declare that they have no conflicts of interest.\newline

\medskip

Data availability statement: no datasets  were generated during and/or analysed during the current study, except to formal computations and one graphic, which it is shown in the paper. But for any doubt or feedback, the corresponding author is available to clarify them on reasonable request.

\end{document}